\documentclass[a4paper,10pt]{article}
\usepackage{amsmath,amssymb}
\usepackage[latin1]{inputenc}
\usepackage{graphicx}

\newcommand{\com}{\operatorname{Comm}}
\newcommand{\perm}{\operatorname{Perm}}
\newcommand{\NAP}{\operatorname{NAP}}
\newcommand{\zero}{\widehat{0}}
\newcommand{\one}{\widehat{1}}

\newcommand{\setA}{\mathtt{A}}
\newcommand{\ttE}{\mathtt{E}}

\newcommand{\compo}[2]{#1 \pmb{(}#2\pmb{)}}

\newcommand{\autgroup}{\operatorname{Aut}}
\newcommand{\aut}{\operatorname{\# Aut}}
\newcommand{\tri}{\triangleleft}

\newcommand{\spec}{\operatorname{Spec}}

\newcommand{\PP}{\mathcal{P}}

\newcommand{\QQ}{\mathbb{Q}}
\newcommand{\FF}{\mathsf{F}}
\newcommand{\GG}{\mathsf{G}}
\newcommand{\sym}{\mathfrak{S}}
\newcommand{\arb}[1]{\includegraphics[height=5mm]{a#1.eps}}
\newcommand{\Farb}[1]{\mathsf{F}_{\big{[}\includegraphics[height=5mm]{a#1.eps}\big{]}}}
\newcommand{\Garb}[1]{\mathsf{G}_{\includegraphics[height=5mm]{a#1.eps}}}
\renewcommand{\tt}[1]{\mathtt{#1}}
\newcommand{\id}{\mathrm{id}}

\newtheorem{theorem}{Theorem}[section] 
\newtheorem{proposition}[theorem]{Proposition} 
 
\newtheorem{corollary}[theorem]{Corollary} 
\newtheorem{lemma}[theorem]{Lemma}

\newenvironment{proof}{\begin{trivlist}\item{\bf{Proof.}}}
  {\hfill\rule{2mm}{2mm}\end{trivlist}}

\title{Relating two Hopf algebras built from an operad}
\author{F. Chapoton and M. Livernet} \date{\today}

\begin{document}

\maketitle

\begin{abstract}
  Starting from an operad, one can build a family of posets. From this
  family of posets, one can define an incidence Hopf algebra. By
  another construction, one can also build a group directly from the
  operad. We then consider its Hopf algebra of functions. 
  We prove that there exists a surjective morphism from the
  latter Hopf algebra to the former one. This is illustrated by the case of
  an operad built on rooted trees, the $\NAP$ operad, where the incidence Hopf algebra is identified
  with the Connes-Kreimer Hopf algebra of rooted trees.
\end{abstract}

\section{Introduction}

Operads were introduced in algebraic topology to deal with loop
spaces, more than 40 years ago. This new algebraic notion has been
somewhat neglected after its introduction, until it appears to be
useful in many other domains, for instance in the algebraic geometry
of moduli spaces of curves, during the 1990's. Since then, there seems
to be a regular activity around operads.

Operads can be defined as objects of any symmetric monoidal category.
Most of the examples considered in the literature live in the category
of sets, of topological spaces, of vector spaces over a field, or of
chain complexes.

The aim of the present article is to relate two different
constructions, both starting from the data of an operad $\PP$ in the
category of sets, of graded commutative Hopf algebras or, from the
dual point of view, of pro-algebraic groups.

The first construction goes as follows. From an operad $\PP$ in the
category of sets, one can define a family of posets, in which the
partial order reflects part of the algebraic structure of the operad.
This partial order has been introduced by Mendez and Yang
\cite{MendezYang91} but rather in the context of species and without
using the term operad; it was rediscovered later by Vallette
\cite{vallette}, who used this to link the Koszul property of the
operad to the Cohen-Macaulay properties of the posets.

Then, one can use this family of posets, which has some adequate
closure property under taking subintervals, as an input to the Schmitt
definition of an incidence Hopf algebra. Therefore, one can build in
this way a first Hopf algebra $H_\PP$ from an operad $\PP$, through
the associated posets.

The second construction of a Hopf algebra from an operad is a direct
one. It is rather the equivalent construction of a pro-algebraic group
$G_\PP$. This has been considered, from different point of views in
\cite{vdlaan,expon,kapr-manin}. As a space, the group $G_\PP$ is an
affine subspace of the completed free $\PP$-algebra on one generator.

Our main general result is the existence of a surjective morphism of
Hopf algebras from the Hopf algebra of functions $\QQ G_\PP$ to the
incidence Hopf algebra $H_\PP$. From the group point of view, this
means that the pro-algebraic group $\spec H_\PP$ is a subgroup of
$G_\PP$.

In the last section of the article, these results are applied to an
operad built on rooted trees, the
operad $\NAP$. We give a precise description of the posets associated
to this operad. We then show that the incidence Hopf algebra for the
$\NAP$ operad is isomorphic to the Hopf algebra of rooted trees which
was introduced by Connes and Kreimer in \cite{ConKre98}. This gives a
second link between this Hopf algebra and operads, after the one
obtained in \cite{prelie} with the pre-Lie operad.

The general theorem is then used, together with a computation of the
Möbius numbers, to find the inverse of a special element of
$G_{\NAP}$. We also provide some other examples of elements of the
group $G_{\NAP}$ and morphisms from this group to more familiar groups
of formal power series in one variable.

\medskip

The present work received support from the ANR grant
BLAN06-1\underline{ }136174.

\section{Set-operads and posets}

Here we recall first the general setting of species and operads, then
the construction of posets starting from an operad (Mendez and Yang
\cite[\S 3.4]{MendezYang91} and Vallette \cite{vallette}) and related results.

\subsection{Species}

The theory of species has been introduced by Joyal \cite{joyal} as a
natural way to deal with generating series. It is closely related to
the notion of $\sym$-module, just as vector spaces are related to
sets.

A \textbf{species} is a functor $\PP$ from the groupoid of finite sets
(the category whose objects are finite sets and morphisms are
bijections) to the category of sets.

For example, the species $\com$ maps a finite set $I$ to the singleton
$\{I\}$ and there is no choice for the bijections.

The category of species is a monoidal category with tensor product
$\circ$ defined by
\begin{equation*}
(F \circ G)(I) =\coprod_{\simeq} F(I/{\simeq}) \times \prod_{J
\in I/\simeq} G(J),
\end{equation*}
where $I$ is a finite set and $\simeq$ runs over the set of
equivalence relations on $I$. Note that this monoidal functor is not
symmetric.

The data of a species $\PP$ is equivalent to the data of a collection of
sets $\PP(n)$ with actions of the symmetric groups. The set $\PP(n)$
can be defined as $\PP(\{1,\dots,n\})$, with the obvious action of the
symmetric group $\sym_n$. The other way round, one can recover the set
$\PP(I)$ as a colimit.

\subsection{Set-operads}

A \textbf{set-operad} $\PP$ is a monoid with unit in the monoidal
category of species for the tensor product $\circ$. This means the
data of a morphism of species
\begin{equation*}
 \gamma : \PP \circ \PP \to \PP,
\end{equation*}
which has to be associative, and a map $e$ from the unit object to $\PP$
satisfying the usual unit axioms. 

An \textbf{augmented operad} $\PP$ is an operad such that
$\PP(\emptyset)$ is empty and the image by $\PP$ of any singleton is a
singleton. We will always assume that the operads we consider are
augmented.

There is an alternative way to describe the composition map $\gamma$
of an operad $\PP$. The data of $\gamma$ as above is equivalent to the
data of maps, for each finite set $I$ and collection of finite sets
$(J_i)_{i\in I}$,
\begin{equation}\label{D-compo}
  \PP(I)\times \prod_{i\in I} \PP(J_i) \to \PP(\coprod_{i\in I} J_i),
\end{equation}
which map $(x,(y_i)_{i \in I})$ to $\compo{x}{(y_i)_{i \in I}}$.

A \textbf{basic set-operad} is a set-operad such that, for each $y\in
\prod_{i\in I} \PP(J_i)$, the map $x \mapsto \compo{x}{y}$ is injective.

\subsection{Posets from set-operads}

Let $\PP$ be a set-operad. Let us denote by $\Pi_\PP$ the species
$\com \circ\, \PP$. Let $I$ be a finite set.

One can build a family of posets on the species $\Pi_\PP$.  More
precisely, there is a partial order on $\Pi_\PP(I)$ for each finite
set $I$ and this construction is functorial in $I$. This means that
the species $\Pi_\PP$ has values in the category of posets rather than
just in the category of sets.

From the definition of $\circ$, one can see that an element $x$ of
$\Pi_{\PP}(I)$ is the data of a partition $\pi_x$ of $I$ and of an
element $x_J$ of $\PP(J)$ for each part $J$ of the partition $\pi_x$.
The definition of the composition maps of $\PP$ in the diagram
(\ref{D-compo}) lifts to the maps
\begin{equation}\label{D-compoext}
\Pi_{\PP}(I)\times \prod_{i\in I} \PP(J_i)\to \Pi_{\PP}(\coprod_{i\in I}
J_i),
\end{equation}
which send $(x=(x_u)_{u\in\pi_x}, y=(y_i)_{i\in I})$ to
\begin{equation*}
  \compo{x}{y}=(\compo{x_u}{(y_i)_{i\in u}})_{u\in\pi_x}.
\end{equation*}
These maps satisfy the following associativity relation
\begin{equation}\label{compo-assoc}
\compo{x}{{\compo{y}{z}}}=\compo{{\compo{x}{y}}}{{z}}.
\end{equation}


Then $x \leq y$ in $\Pi_\PP(I)$ if there exists an element
$\theta(x,y) \in \Pi_{\PP}(\pi_x)$ such that
\begin{equation*}
  \compo{\theta(x,y)}{ x } = y.
\end{equation*}
Note that this definition implies that the partition $\pi_x$ is finer
than the partition $\pi_y$.

The poset $\Pi_{\PP}(I)$ has a unique minimal element, 
denoted by $\zero$. 

The following proposition is statement 3. in \cite[Thm. 3.4]{MendezYang91}.

\begin{proposition}\label{L-orderpreserving} 
  Let $\PP$ be a basic set-operad. Let $x\in\Pi_{\PP}(I)$. The poset
  $\{y\in\Pi_{\PP}(I)|\ x\leq y\}$ is isomorphic to the poset
  $\Pi_{\PP}(\pi_x)$.
\end{proposition}

\begin{proof} 
  Since $\PP$ is a basic set-operad, if $x\leq y$ there is a unique
  $\theta(x,y)\in \Pi_{\PP}(\pi_x)$ such that
  $\compo{\theta(x,y)}{x}=y$.  The bijection sends $y$ to
  $\theta(x,y)$. The inverse map sends $a\in \Pi_{\PP}(\pi_x)$ to
  $\compo{a}{x}$. Assume $x\leq y\leq z$. By definition
  $\compo{\theta(x,y)}{x}=y$ and $\compo{\theta(x,z)}{x}=z$. Since
  $y\leq z$ there is a unique $\theta(y,z)\in\Pi_{\PP}(\pi_y)$ such
  that $\compo{\theta(y,z)}{y}=z$.  As a consequence, using the
  associativity relation (\ref{compo-assoc}), one has
\begin{equation*}
  z=\compo{\theta(y,z)}{y}=\compo{\theta(y,z)}{\compo{\theta(x,y)}{x}}=\compo{\compo{\theta(y,z)}{\theta(x,y)}}{x}
=\compo{\theta(x,z)}{x}.
\end{equation*}
The uniqueness of the elements $\theta$ permits to conclude that
$\compo{\theta(y,z)}{\theta(x,y)}=\theta(x,z)$
and
$\theta(x,y)\leq\theta(x,z).$

Conversely, if $a\leq b\in \Pi_{\PP}(\pi_x)$ one has clearly 
$\compo{a}{x}\leq \compo{b}{x}$.
\end{proof}

\medskip

If this construction is applied to the set-operad $\com$, the poset
$\Pi_{\com}(I)$ is the usual partial order by refinement on the
partitions of the set $I$.

One can similarly get the poset of pointed-partitions, when this
construction is applied to the set-operad $\perm$ \cite{chapoval}.

\medskip

Vallette has used these posets to give a Koszulness criteria for
operads. Let us just recall the result here.

\begin{proposition}[\cite{vallette}, Theorem 12]
  \label{val_critere}
  Let $\PP$ be a set-operad which is basic, quadratic and augmented.
  Then the associated linear operad $\QQ\PP$ is Koszul if and only if
  all maximal intervals in $\Pi_{\PP}(I)$ are Cohen-Macaulay for all
  $I$.
\end{proposition}

\section{Incidence Hopf algebras}

Here we recall briefly a construction of William Schmitt
\cite{schmitt} building a commutative Hopf algebra from a family of
posets satisfying some conditions. We then derive our first Hopf
algebra built from an operad from the composition of the construction
of Mendez-Yang and Vallette with this construction of Schmitt.

\subsection{Good families of posets}

Suppose we are given a collection of posets $(P_\alpha)_{\alpha \in
  \setA}$. The collection $(P_\alpha)_{\alpha\in\setA}$ is called a
\textbf{good collection} if it satisfies the following conditions.
\begin{enumerate}
\item Each poset $P_\alpha$ has a minimal element $\zero$ and a maximal
  element $\one$ (it is an \textbf{interval}).
\item For all $\alpha \in \setA$ and all $x$ in $P_\alpha$, the interval
  $[\zero,x]$ is isomorphic to a product of posets $\prod_\beta P_\beta$ and
  the interval $[x,\one]$ is isomorphic to a product of posets
  $\prod_\gamma P_\gamma$.
\end{enumerate}

As a simple example of good collection, one can consider the family
of all total orders. Another example is the family of boolean posets.

Remark: it follows from this definition applied to the interval
$[\zero,\zero]$ in any poset $P_\alpha$ that a good collection
contains at least one poset $P_{\varepsilon}$ with only one element.

\subsection{Hopf algebra from a good collection}

Let $(P_\alpha)_{\alpha \in\setA}$ be a good collection. Let us
consider the collection of all finite products $\prod_\beta P_\beta$.
Let us denote by $\bar{\setA}$ this larger set of posets.

The collection $\bar{\setA}$ of posets is closed under products by
construction. It is also closed under taking initial intervals
$[\zero,x]$ or final intervals $[x,\one]$. Hence it is also closed
under taking any subinterval, because any interval $[x,y]$ is a final
interval in the initial interval $[\zero,y]$.

A collection of posets which is closed under products and closed under
taking subintervals is called a \textbf{hereditary collection} in
\cite{schmitt}. The collection $\bar{\setA}$ is therefore a hereditary
collection.

Let us denote by $[\setA]$ the set of isomorphism classes of posets in
$\setA$, and by $[\bar{\setA}]$ the set of isomorphism classes of
posets in $\bar{\setA}$. Elements in these sets will be denoted by
$[\alpha],[\beta],\dots$, which will also mean the isomorphism class
of $\alpha,\beta,\dots \in \setA$ or $\bar{\setA}$.

One can then consider the vector space $H_\setA$ with a basis
$\FF_{[\alpha]}$ indexed by the set $[\bar{\setA}]$.

Then $H_\setA$ is a commutative algebra for the product induced by the
direct product of posets:
\begin{equation*}
  \FF_{[\alpha]}\FF_{[\beta]}=\FF_{[\alpha\times \beta]}.
\end{equation*}
This algebra is generated by the elements $\FF_{[\alpha]}$ with
$[\alpha] \in [\setA]$. Note that one can remove the unit
$\FF_{[\varepsilon]}$ from this set of generators. The algebra $H_\setA$
may not be free on this reduced set of generators, as there can be
isomorphisms $\prod_\beta P_\beta \simeq \prod_\gamma P_\gamma$ with a
different number of factors or with non pairwise-isomorphic factor
posets.

The space $H_\setA$ is also a coalgebra for a coproduct $\Delta$ whose
value on the generator $ \FF_{[\alpha]}$ is
\begin{equation*}
  \Delta(\FF_{[\alpha]}) = \sum_{x \in P_\alpha} \FF_{[\zero,x]} 
  \otimes \FF_{[x,\one]},
\end{equation*}
where the intervals in indices stand for their isomorphism classes.

In fact, this formula is enough to define the coproduct $\Delta$, which is
compatible with the product on $H_\setA$.

To summarize,
\begin{proposition}
  The space $H_{\setA}$ endowed with its commutative product and the
  coproduct $\Delta$ is a commutative Hopf algebra. The unit is
  $\FF_{[\varepsilon]}$, where $[\epsilon]$ is the isomorphism class of
  the singleton interval.
\end{proposition}

This is a consequence of the general theorem of Schmitt on hereditary
collections of posets \cite[Theorem 4.1]{schmitt}.

\subsection{Group from a good family}

The commutative Hopf algebra $H_{\setA}$ is the space of functions on
a pro-algebraic group $\spec H_{\setA}$, the elements of which can be
seen as some formal power series indexed by elements of $[\setA]$. The
fact that $H_{\setA}$ is not necessarily a polynomial algebra on the
set $(\FF_{[\alpha]})_{[\alpha]\in[\setA]}$ is equivalent to the
possible existence of some universal relations between the
coefficients of these series (see Lemma \ref{decrire_spec} for an
instance of this phenomenon). The fact that $\FF_{[\varepsilon]}$ is
the unit means that the coefficient of $[\varepsilon]$ in these series
is $1$.

An element of this pro-algebraic group can be considered as a function
on the collection of isomorphism classes of posets
$(P_{[\alpha]})_{[\alpha]\in[\setA]}$. The product in the group
provides information on the posets, by the classical theory of
M{\"o}bius functions, zeta functions and incidence algebra of posets,
see \cite{rota,stanley1}.

The following proposition gives an example of computation in the
pro-algebraic group $\spec H_{\setA}$.

\begin{proposition}[\cite{schmitt} \S 7]
  \label{zeta_inverse} 
  The group product in $\spec H_{\setA}$ gives the usual convolution
  product on functions over the posets $P_{[\alpha]}$ for $[\alpha] \in[\setA]$. Consider
  in $\spec H_{\setA}$ the Möbius series
  \begin{equation*}
    \tt{M}=\sum_{[\alpha]\in [\setA]} \mu(P_{[\alpha]}) [\alpha],
  \end{equation*}
  where $\mu(P_{[\alpha]})$ is the Möbius number of the poset $P_{[\alpha]}$, and the
  Zeta series
  \begin{equation*}
    \tt{Z}=\sum_{[\alpha]\in [\setA]} [\alpha].  
  \end{equation*}
  Then $\tt{M}$ is the inverse of $\tt{Z}$ in $\spec H_{\setA}$.
\end{proposition}


\subsection{From operads to incidence Hopf algebras}
\label{incidence_from_operad}

Here we show that one can use the posets $\Pi_\PP(I)$ associated
with a basic set-operad $\PP$ to define an incidence Hopf algebra
$H_\PP$ by using Schmitt construction for an hereditary family.

Indeed the intervals in $\Pi_\PP(I)$ are products of minimal intervals
as stated in the following proposition.

\begin{proposition}
  \label{sous-intervalles}
 Let $y \in \PP(n)$ for some $n$. Let $\zero \leq x \leq y$. Assume
that $x$ has components $(x_u)_{u \in \pi_x}$. 
The interval $[\zero,x]$ is isomorphic to the product of posets
  $\prod_{u \in \pi_x} [\zero,x_u]$. The interval $[x,\one]=[x,y]$ is
  isomorphic to the poset $ [\zero,\theta(x,y)] $, where $\theta(x,y)$
  is the unique element of $\PP(\pi_x)$ such that
  $\compo{\theta(x,y)}{x}=y$.
\end{proposition}

\begin{proof}
  The isomorphism between $[\zero,x]$ and $\prod_{u \in    \pi_x}
  [\zero,x_u]$ is a direct consequence of the definition of
  the partial order. Indeed, one has $z \leq x$ if and only if the
  partition $\pi_z$ is finer than the partition $\pi_x$ and for each
  part $u$ of $\pi_x$, one has $z_u \leq x_u$, where $z_u$ and $x_u$
  are the restrictions of $z$ and $x$ to $u$. This allows to prove the
  expected isomorphism.

  Let us now consider the interval $[x,y]$. The order preserving
  isomorphism of Prop. \ref{L-orderpreserving} between $\{z|
  x\leq z\}$ and $\Pi_\PP(\pi_x)$ induces an isomorphism between the
  intervals $[x,y]$ and $[\zero,\theta(x,y)]$.
\end{proof}

As a consequence, if $x\leq y \in \Pi_\PP(I)$ and
$y=(y_u)_{u\in\pi_y}$ then
\begin{equation}\label{intervalles}
[x,y] \simeq \prod_{u\in\pi_y} [\zero, \theta(x_u,y_u)].
\end{equation}

Let $\setA_\PP$ be the set of coinvariants for the species $\PP$. 
For each coinvariant $\alpha\in \PP(n)_{\sym n}$, let $r(\alpha)$ be a representative of
$\alpha$ in $\PP(n)$. Let us define a poset $P_{\alpha}$ as the
interval $[\zero,r(\alpha)]$ in $\Pi_\PP(n)$.

\begin{proposition}
  The collection of posets $(P_\alpha)_{\alpha\in{\setA_\PP}}$ is a good
  family of posets. The resulting incidence Hopf algebra is denoted by $H_\PP$.
\end{proposition}

\begin{proof}
Obviously, all posets $P_\alpha$ are intervals. There remains only to
prove the stability property. For any $\alpha\in\setA_\PP$, two
representatives $r(\alpha)$ and $s(\alpha)$ give two isomorphic
intervals  $[\zero,r(\alpha)]$ and $[\zero,s(\alpha)]$ in
$\Pi_\PP(n)$. Thus, the stability property follows from Prop. \ref{sous-intervalles}.
\end{proof}

\section{Groups from operads}

\label{group_from_operad}

Here we recall the construction of a group from an operad. The Hopf
algebra of its functions gives
our second Hopf algebra built from an operad. 

We will work with a set-operad $\PP$, but the construction is just the same
for an operad in the category of vector spaces. This simple
construction has already been considered from different viewpoints in
\cite[Chap. I, \S 1.2]{vdlaan} and \cite{expon,kapr-manin}.

Let $\PP$ be an augmented set operad. In this section, we will use the
description of a species $\PP$ as a collection of modules $\PP(n)$
over the symmetric groups.

Let $\QQ\setA_{\PP}=\oplus_n
\QQ\PP(n)_{\sym_n}$ be the direct sum of the coinvariant spaces, which
can be identified with the underlying vector space of the free
$\PP$-algebra on one generator, and $\widehat{\QQ\setA}_{\PP}={\prod}_n
\QQ\PP(n)_{\sym_n}$ be its completion.

Let $\alpha=\sum_m \alpha_m$, $\beta=\sum_n \beta_n$ be two elements
of $\widehat{\QQ\setA}_{\PP}$ with
$\alpha_m,\beta_m$ elements of $\QQ\PP(m)_{\sym_m}$. 
Choose any representatives
${x}_m=r(\alpha_m)$ of $\alpha_m$ (resp.  ${y}_m=r(\beta_m)$ of
$\beta_m$) in $\QQ\PP(m)$. Then one can check that the following formula defines a
product on $\widehat{\QQ\setA}_{\PP}$:
\begin{equation}
  \label{product}
  \alpha \times \beta = \sum_{m\geq 1} \, \sum_{n_1,\dots,n_m \geq 1} 
  \langle \compo{{x}_m} {{y}_{n_1},\dots,{y}_{n_m} }\rangle,
\end{equation}
where $\langle\,\rangle$ is the quotient map to the coinvariants and
$(x,y_1,\dots,y_m)\mapsto \compo{x}{y_1,\dots,y_m}$ is the composition
map of the operad $\PP$.

\begin{proposition}
  The product $\times$ defines the structure of an associative monoid
  on the vector space $\widehat{\QQ\setA}_{\PP}$. Furthermore, this product is
  $\QQ$-linear on its left argument.
\end{proposition}
\begin{proof}
  Let us first prove the associativity. Let $\delta=\sum_p \delta_p$
  and fix representatives $z_p=r(\delta_p)$. On the one hand, one has
  \begin{multline}
    (\alpha\times \beta)\times \delta=\sum_m \sum_{p_1,\dots,p_m}
    \langle\compo{r((\alpha\times \beta)_m)}{ {z}_{p_1}, \dots,
    {z}_{p_m} }\rangle\\=
    \sum_m\sum_{n_1,\dots,n_m}\,\sum_{p_1,\dots,p_{n_1+\dots+n_m}}
    \langle\compo{\compo{{x}_m}{ {y}_{n_1},\dots,
    {y}_{n_m} }} {{z}_{p_1},
    \dots, {z}_{p_{n_1+\dots+n_m}} }\rangle.
  \end{multline}

  On the other hand, one has 
  \begin{multline}
    \alpha\times (\beta \times
    \delta)=\sum_m\sum_{n_1,\dots,n_m}\langle\compo{{x}_m} {r((\beta \times
      \delta)_{n_1}),\dots,r((\beta\times \delta)_{n_m}) }\rangle
    \\=\sum_m\sum_{n_1,\dots,n_m}\sum_{(q_{i,j})}\langle\compo{{x}_m}{
    \compo{{y}_{n_1}}{
    {z}_{q_{1,1}},\dots,{z}_{q_{1,n_1}}},\dots,\compo{{y}_{n_m}} {{z}_{q_{m,1}},\dots,{z}_{q_{m,n_m}}}}\rangle.
  \end{multline}
  Using then the ``associativity'' of the operad, one gets the
  associativity of $\times$. It is easy to check that the image
  $\varepsilon$ of the unit $e$ of the operad $\PP$ is a two-sided
  unit for the $\times$ product. The left $\QQ$-linearity is clear
  from the formula (\ref{product}).
\end{proof}

\begin{proposition}
  An element $\beta$ of $\widehat{\QQ\setA}_{\PP}$ is invertible for $\times$ if
  and only if the first component $\beta_1$ of $\beta$ is non-zero.
\end{proposition}
\begin{proof}
  The direct implication is trivial. The converse is proved by a
  very standard recursive argument.
\end{proof}

Let us call $G_{\PP}$ the set of elements of $\widehat{\QQ\setA}_{\PP}$ whose
first component is exactly the unit $\varepsilon$. This is a subgroup
for the $\times$ product of the set of invertible elements.

\begin{proposition}
  The construction $G$ is a functor from the category of augmented
  operads to the category of groups.
\end{proposition}
\begin{proof}
  The functoriality follows from inspection of the definitions of
  $\widehat{\QQ\setA}_{\PP}$ and $\times$.
\end{proof}

In fact, one can see $G_\PP$ as the group of $\QQ$-points of a
pro-algebraic and pro-unipotent group. The Lie algebra of this
pro-algebraic group is given by the usual linearization process on the
tangent space (an affine subspace of $\widehat{\QQ\setA}_{\PP}$),
resulting in the formula
\begin{equation*}
  [\alpha, \beta] = \sum_{m\geq 1} \sum_{n \geq 1} 
  \langle {x}_{m} \circ {y}_{n}
    -{y}_{n} \circ {x}_{m}\rangle,
\end{equation*}
where
\begin{equation*}
  {x}_m \circ {y}_n=\sum_{i=1}^{m} \compo{{x}_m}{
\underbrace{e,\dots,e}_{i-1\text{ units}},{y}_n,e,\dots,e}.
\end{equation*}

The graded Lie algebra structure on $\QQ\setA_{\PP}$ defined by the same
formulas has already appeared in the work of Kapranov and Manin on the
category of right modules over an operad \cite[Th.
1.7.3]{kapr-manin}.

\medskip


The Hopf algebra $\QQ[G_\PP]$ of functions on $G_\PP$ is the free
commutative algebra generated by $\GG_\alpha$ for $\alpha$ in the set
$\setA_\PP$ but the unit invariant $\varepsilon$. An element $g$ of
$G_\PP$ can be seen as a formal sum
\begin{equation*}
  g=\sum_{\alpha \in \setA_\PP} \GG_{\alpha}(g) \alpha,
\end{equation*}
where $\GG_{\varepsilon}=1$.
As a function on $G_{\PP}$, the value of $\GG_\alpha$ on an element
$g$ of $G_{\PP}$ is the coefficient of $\alpha$ in the expansion of
$g$.

\section{Main theorem}
\label{S-maintheorem}

Here we show that the incidence Hopf algebra $H_{\PP}$ defined in \S
\ref{incidence_from_operad} is a quotient
of the Hopf algebra of functions $\QQ[G_\PP]$ on the group of formal
power series defined directly from the operad $\PP$ by the
construction of \S \ref{group_from_operad}.

This also means that the group $\spec H_{\PP}$ is a subgroup of the
group $G_\PP$.

Let us consider the coproduct $\Delta$ in the incidence Hopf algebra
$H_\PP$. This space has a basis indexed by the set $[\bar{\setA}_\PP]$
of isomorphism classes of products of posets. The set $[\setA_\PP]$ is
a subset of $[\bar{\setA}_\PP]$. If one considers the coproduct on one
element $\FF_{[\alpha]}$ with $\alpha \in [\setA_\PP]$, then it can be
written uniquely as a linear combination
\begin{equation*}
  \Delta(\FF_{[\alpha]})=\sum_{[\gamma],[\beta]}
  \mathbf{f}_{[\alpha]}^{[\beta],[\gamma]} \FF_{[\beta]} \otimes \FF_{[\gamma]},
\end{equation*}
where $([\alpha],[\beta],[\gamma])$ in $[\setA_\PP] \times
[\bar{\setA}_\PP] \times [\setA_\PP]$. Indeed, the fact that this sum
only runs over $\gamma\in [\setA_\PP]$ (and not $[\bar{\setA}_\PP]$)
follows from the description of the subintervals in Prop.
\ref{sous-intervalles}.

Therefore, for each triple $([\alpha],[\beta],[\gamma])$ in $[\setA_\PP] \times
[\bar{\setA}_\PP] \times [\setA_\PP]$, one can define a coefficient
$\mathbf{f}_{[\alpha]}^{[\beta],[\gamma]}$ by the previous expansion.

Similarly, one can consider the Hopf algebra of functions on the group
$G_\PP$ and define, for each triple $(\alpha,\beta,\gamma)$ with
$\alpha$ an element of $\PP(n)_{\sym_n}$ for some $n$, $\gamma$ an
element of $\PP(k)_{\sym_k}$ for some $k\leq n$ and $\beta$ an element
of $(\Pi_{ \PP})(n)_{\sym_n}$ with $k$ parts, a coefficient
$\mathbf{g}_{\alpha}^{\beta,\gamma}$ by
\begin{equation*}
  \Delta(\GG_\alpha)=\sum_{\gamma,\beta}
  \mathbf{g}_{\alpha}^{\beta,\gamma} \GG_\beta \otimes \GG_\gamma,
\end{equation*}
where $\GG_\beta$ is the product $\prod_t \GG_{\beta_t}$
over the set of components of $\beta$.

\medskip

Let us choose for the rest of this section a triple
$(\alpha,\beta,\gamma)$ as above. We will compare the coefficients
$\mathbf{f}_{[\alpha]}^{[\beta],[\gamma]}$ and $\mathbf{g}_{\alpha}^{\beta,\gamma}$.

Let us denote by $\langle\,\rangle$ the projections to coinvariants from
$\PP(n)$ to $\PP(n)_{\sym_n}$.

Let us pick a representative $r(\alpha)$ of $\alpha$ in $\PP(n)$ and a
representative $r(\gamma)$ of $\gamma$ in $\PP(k)$. Let us also choose
a representative $r(\beta)$ of $\beta$ in $(S^k \PP)(n)$ with the
following property: the partition of $\{1,\dots,n\}$ induced by the
components of the representative $r(\beta)$ is the standard partition
\begin{equation*}
  p_{\operatorname{std}}=\{1,\dots,\ell_1\}\sqcup \{ \ell_1+1,\dots,\ell_1+\ell_2\}\sqcup \dots
  \sqcup \{\ell_{1}+\dots+\ell_{k-1}+1,\dots,\ell_1+\dots+\ell_k \}.
\end{equation*}

This allows to define a bijection between the set of components of
$\beta$ and the set $\{1,\dots,k\}$. Then one will denote by $\beta_i$
the component indexed by $i$. By the unique increasing renumbering,
this also gives representatives $r(\beta_i)$ of $\beta_i$ in
$\PP(\ell_i)$.

Let us introduce the automorphism groups $\autgroup(\alpha)$,
$\autgroup(\gamma)$ and $\autgroup(\beta)$. They are rather the
automorphisms groups of representatives $r(\alpha)$, $r(\beta)$ and
$r(\gamma)$. The group $\autgroup(\beta)$ decomposes into a
semi-direct product
\begin{equation*}
  \autgroup(\beta) = \left(\prod_{i=1}^{k} \autgroup(\beta_i)\right) \rtimes \autgroup_0(\beta),
\end{equation*}
where $\autgroup_0(\beta)$ is a subgroup of the permutation group $\sym_k$
of the set of components of $\beta$.

From the description of the coproduct in the incidence Hopf algebra,
the coefficient $\mathbf{f}_{[\alpha]}^{[\beta],[\gamma]}$ is the
cardinal number of the following set
\begin{equation}
  \label{cardf}
  \{p, \sigma\in \sym_k ,u,v_i \mid r(\alpha)=\compo{u}{
  v_{1},\dots,v_{k}},\quad \langle u \rangle=\gamma,\quad \langle
  v_{\sigma(i)} \rangle=\beta_i\},
\end{equation}
where $p$ is a partition of $\{1,\dots,n\}$ with $k$ parts $p_i$
ordered by their least element, $\sigma\in \sym_k$, $u\in \PP(k)$ and
$v_i\in \PP(p_i)$ for $i=1,\dots,k$.

Let us introduce the set $\ttE_\mathbf{f}(\alpha,\beta,\gamma)$
consisting of
\begin{equation*}
  \{p,\sigma,\psi,\phi_i,u, v_i \mid r(\alpha)=\compo{u}
  {v_{1},\dots,v_{k}},\quad u \simeq^{\psi} r(\gamma),\quad v_{\sigma(i)}
  \simeq^{\phi_i} r(\beta_i) \}
\end{equation*}
where $p$ is a partition of $\{1,\dots,n\}$ with $k$ parts $p_i$
ordered by their least element, $\sigma\in \sym_k$, $u\in \PP(k)$,
$v_i\in \PP(p_i)$ for $i=1,\dots,k$, $\psi\in\sym_k$ and $\phi_i$ is
bijection from the part $p_{\sigma(i)}$ to the set $\{1,\dots,\ell_i\}$.

\begin{proposition}
  \label{decrire_ensf}
  The set $\ttE_\mathbf{f}(\alpha,\beta,\gamma)$ satisfies
  \begin{equation*}
    \aut(\beta) \aut(\gamma) \,\mathbf{f}_{[\alpha]}^{[\beta],[\gamma]}= \#\ttE_\mathbf{f}(\alpha,\beta,\gamma).
  \end{equation*}
\end{proposition}
\begin{proof}
  The group $\autgroup(\beta)\times \autgroup(\gamma)$ acts freely on
  $\ttE_\mathbf{f}(\alpha,\beta,\gamma)$ and the orbits are in
  bijection with the set described in (\ref{cardf}) whose cardinality is
  $\mathbf{f}_{[\alpha]}^{[\beta],[\gamma]}$.
\end{proof}

\medskip

From the description of the product in the group $G_{\PP}$, the
coefficient $\mathbf{g}_{\alpha}^{\beta,\gamma}$ is the cardinal of
the following set
\begin{equation}
  \label{cardg}
  \{ \tau\in\sym_k/\autgroup_0(\beta) \mid \alpha=\langle
  \compo{r(\gamma)}{r(\beta_{\tau(1)}),\dots,r(\beta_{\tau(k)})}
  \rangle \}.
\end{equation}

Let us introduce the set $\ttE_\mathbf{g}(\alpha,\beta,\gamma)$ consisting of 
\begin{equation*}
  \{ \tau\in\sym_k, \phi \in \sym_n \mid r(\alpha) \simeq^{\phi} \compo{r(\gamma)}{r(\beta_{\tau(1)}),\dots,r(\beta_{\tau(k)})} \}).
\end{equation*}
 
\begin{proposition}
  \label{decrire_ensg}
  The set $\ttE_\mathbf{g}(\alpha,\beta,\gamma)$ satisfies
  \begin{equation*}
    \aut(\alpha) \aut_0(\beta) \,\mathbf{g}_{\alpha}^{\beta,\gamma} = \#\ttE_\mathbf{g}(\alpha,\beta,\gamma).
  \end{equation*}  
\end{proposition}
\begin{proof}
  The group $\autgroup(\alpha)\times \autgroup_0(\beta)$ acts freely on
  $\ttE_\mathbf{g}(\alpha,\beta,\gamma)$ and the orbits are in
  bijection with the set described in (\ref{cardg}) whose cardinality is
  $\mathbf{g}_{\alpha}^{\beta,\gamma}$.
\end{proof}

\medskip

Let us now show that the sets $\ttE_\mathbf{f}$ and $\ttE_\mathbf{g}$ are
just the same.

\begin{proposition}
  \label{ensf_ensg}
  There is a bijection between $\ttE_\mathbf{f}(\alpha,\beta,\gamma)$
  and $\ttE_\mathbf{g}(\alpha,\beta,\gamma)$.
\end{proposition}

\begin{proof}
  Recall the definition of the set
  $\ttE_\mathbf{f}(\alpha,\beta,\gamma)$ consisting of
  \begin{equation*}
    \{ p,\sigma,\psi,\phi_i, u, v_i \mid r(\alpha)=\compo{u}{v_{1},\dots,v_{k}},\quad u \simeq^{\psi} r(\gamma),\quad v_{\sigma(i)}
    \simeq^{\phi_i} r(\beta_i) \}.
  \end{equation*}

  Let us pick an element in this set. Then there exists a unique
  permutation $\phi\in\sym_n$ induced by the collection of bijections
  $\phi_i$. This bijection maps the partition $p$ to the standard
  partition $p_{\operatorname{std}}$, changing the order of the parts
  according to $\sigma$. It provides an isomorphism between $r(\alpha)$ and 
  \begin{equation*}
    \compo{\sigma^{-1}(u)}{r(\beta_1),\dots,r(\beta_k)}.
  \end{equation*}
  Then one can use $\psi$ and $\sigma^{-1}$ to define a unique isomorphism $\tau$
  between $\sigma^{-1}(u)$ and $ r(\gamma)$.
  This gives us an equality 
  \begin{equation*}
    \compo{\sigma^{-1}(u)}{r(\beta_1),\dots,r(\beta_k)}=  \compo{  r(\gamma)}{ r(\beta_{\tau(1)}),\dots,r(\beta_{\tau(k)})},
  \end{equation*}
hence a unique element in the set
    \begin{equation*}
  \ttE_\mathbf{g}(\alpha,\beta,\gamma)=   \{ \tau\in\sym_k, \phi\in
  \sym_n \mid r(\alpha) \simeq^{\phi} 
\compo{r(\gamma)}{r(\beta_{\tau(1)}),\dots,r(\beta_{\tau(k)}} \}.
  \end{equation*}
\end{proof}

One can now prove the existence of a morphism between Hopf algebras or
equivalently between groups.

\begin{theorem}
  \label{main_thm}
  The map $\rho : \GG_\alpha \mapsto
  \frac{\FF_{[\alpha]}}{\aut(\alpha)}$ defines a surjective
  morphism from the Hopf algebra $\QQ[G_\PP]$ of coordinates on the
  group $G_\PP$ to the incidence Hopf algebra $H_\PP$. In terms of
  groups, this means that the group $\spec H_{\PP}$ is a subgroup of
  the group $G_{\PP}$.
\end{theorem}

\begin{proof}
  The Hopf algebra $\QQ[G_\PP]$ is commutative and freely generated by
  the set of coinvariants of $\PP$ (but the unit). On the other hand,
  the incidence Hopf algebra is commutative and generated by the
  isomorphism classes of maximal intervals (but the trivial interval).

  As intervals coming from the same coinvariant are obviously
  isomorphic, the proposed map is well defined from the set of
  coinvariants to the set of isomorphism classes of intervals. Then one
  can uniquely extend this map into a morphism of algebras, because
  $\QQ[G_\PP]$ is a free commutative algebra. This morphism is
  surjective by construction.

  According to the notation introduced before, we have to prove that
  \begin{equation*}
    \aut(\alpha) \,\mathbf{g}_{\alpha}^{\beta,\gamma} =   \prod_t \aut(\beta_t) \aut(\gamma) \,\mathbf{f}_{[\alpha]}^{[\beta],[\gamma]}.
  \end{equation*}

  By the semi-direct product structure of $\autgroup(\beta)$, this is
  equivalent to
  \begin{equation*}
    \aut(\alpha) \aut_0(\beta) \,
    \mathbf{g}_{\alpha}^{\beta,\gamma} =  \aut(\beta) \aut(\gamma) \,\mathbf{f}_{[\alpha]}^{[\beta],[\gamma]}.    
  \end{equation*}

  This follows in turn from Prop. \ref{decrire_ensf},
  Prop. \ref{decrire_ensg} and Prop. \ref{ensf_ensg}.
\end{proof}

Let us consider briefly a simple example, which is the operad $\com$.
For each $n$, the space $\com(n)$ is the trivial module over the
symmetric group $\sym_n$, hence $\com(n)_{\sym_n}$ has dimension $1$.
The algebra of functions $\QQ[G_{\com}]$ is free on one generator
$\GG_n$ in each degree $\geq 2$. Using the definition of $G_{\com}$,
one can check that the group $G_{\com}$ is isomorphic to the group of
formal power series
\begin{equation*}
  f = x+\sum_{n \geq 2} \GG_n(f) x^n
\end{equation*}
for composition (a group of formal diffeomorphisms).

On the other hand, there is only one interval in $\Pi_{\com}(n)$,
which is the usual poset of partitions. The incidence Hopf algebra of
this family of posets is very classical \cite[Ex. 14.1]{schmitt},
freely generated by one element $\FF_n$ in each degree and isomorphic
to the Faà di Bruno Hopf algebra, which is the Hopf algebra of
functions on the group of formal power series
\begin{equation*}
  f = x+\sum_{n \geq 2} \FF_n(f) \frac{x^n}{n!}.
\end{equation*}
for composition.

Hence, in the case of $\com$, the morphism from $\QQ[G_{\com}]$ to
$H_{\com}$ which maps $\GG_n$ to $\FF_n /n!$ is an isomorphism. The
next section is devoted to the case of the operad $\NAP$ where the
surjective morphism is not an isomorphism.

\section{Application to the $\NAP$ operad}

\subsection{The $\NAP$ operad}

The us first recall the definition of the $\NAP$ operad, which has
been introduced in \cite{rigidity}. The name $\NAP$ stands
for "non-associative permutative".

Let $I$ be a finite set. The set $\NAP(I)$ is the set of rooted trees
with vertices $I$, that is, connected and simply connected graphs with
a distinguished vertex called the root. The unit is the unique rooted
tree on the set $\{i\}$ for any singleton. 

We use the notation
\begin{equation*}
  t=B(r,t_1,\dots,t_k)
\end{equation*}
for a rooted tree $t$ built from the rooted trees $t_i$ by adding an edge from the
root of each rooted tree $t_i$ to a disjoint vertex $r$, which becomes the
root of $t$.

Let us describe the composition $\compo{t}{(s_i)_{i \in I}}$, where
$t \in \NAP(I)$ and $s_i \in \NAP(J_i)$.

Consider the disjoint union of the rooted trees $s_i$ and add some edges: for
each edge of $t$ between $i$ and $i'$ in $I$, add an edge between the
root of $s_i$ and the root of $s_{i'}$. The result is a rooted tree on the
vertices $\sqcup_i J_i$. This is $\compo{t}{(s_i)_{i \in I}}$. The
root of this rooted tree is the root of $s_k$ where $k$ is the index of the
root of $t$.

A $\NAP$-algebra is a vector space $V$ endowed with a bilinear map
$\tri$ from $V\otimes V \to V$ such that
\begin{equation*}
  (a \tri b) \tri c = (a \tri c) \tri b.
\end{equation*}

The free $\NAP$-algebra on a set $\tt{S}$ of generators has a basis
indexed by rooted trees together with a bijection from vertices to
$\tt{S}$.  The product $s \tri t$ of two such rooted trees is obtained by
grafting the root of $t$ on the root of $s$: $B(r,s_1,\ldots,s_k)\tri
t=B(r,s_1,\ldots,s_k,t)$.

Let us note here that $\NAP$ is a basic set-operad. Indeed, one can
recover $t$ from $u=\compo{t}{(s_i)_{i\in I}}$ and the collection
$(s_i)_{i\in I}$ as the restriction of $u$ to the vertices which are
roots of some rooted tree $s_i$.

\subsection{Posets associated with $\NAP$}

Let us describe the posets $\Pi_{\NAP}(I)$. The set $\Pi_{\NAP}(I)$
consists of forests of rooted trees with vertices labeled by $I$.

The covering relations can be described as follows: a forest $x$ is
covered by a forest $y$ if $y$ is obtained from $x$ by grafting the
root of one component of $x$ on the root of another component of $x$.
In the other direction, $x$ is obtained from $y$ by removing an edge
incident to the root of one component of $y$.

By relation (\ref{intervalles}), any interval in $\Pi_{\NAP}(I)$ is a
product of intervals of the form $[\zero, t]$ for $t\in\NAP(J)$.

Let us introduce the following order on rooted trees: $t\leq_s t'$ or
$t'$ is a \textbf{sub-rooted tree} of $t$ if $t'$ is the restriction
of $t$ to a subset of vertices containing the root of $t$, such that
every vertex of $t$ lying on the path between the root and a
vertex of $t'$ is also in $t'$. If $t$ itself is seen as a poset with
its root as minimum element, then this just means that $t'$ is a lower
ideal of $t$.

Let $\one$ denotes the root of $t$ and $[t,\one]_s$ be the interval
between $t$ and $\one$ for the order $\leq_s$.  A rooted tree $t$ is covered
by a rooted tree $t'$ if $t'$ is obtained from $t$ by removing a leaf.

\begin{proposition}
  \label{subtree}
  The interval $[\zero,t]$ is isomorphic to the interval $[t,\one]_s$.
\end{proposition}

\begin{proof}
  Let $x$ be a forest such that $x\leq t$ and let $r_x$ be the family
  of its roots. By definition of the order relation, there exists a
  rooted tree $z$ such that $t=\compo{z}{x}$. It means that the vertices
  indexed by $r_x$ form a sub-rooted tree of $t$.  The isomorphism
  from $[\zero,t]$ to $[t,\one]_s$ sends $x$ to this sub-rooted tree.
  The inverse morphism is the following: let $t'$ be a sub-rooted tree
  of $t$ and $r_{t'}$ the set of its vertices. Again, in view of the
  composition, there is a unique forest $x\leq t$ whose roots are
  indexed by $r_{t'}$ such that $t=\compo{t'}{x}$.
\end{proof}

This proposition in fact shows that $[\zero,t]$ is isomorphic to the
lattice of lower order ideals of the rooted tree $t$ seen as a poset. From
this, it follows that $[\zero,t]$ is a distributive lattice
(\cite{birkhoff}). According to \cite[Example 2.4]{mcnamara}), this
also implies Prop. \ref{supersolvable} below.


\begin{proposition}
  \label{supersolvable}
  The intervals in $\Pi_{\NAP}(n)$ are $\sym_n$ EL-shellable and
  supersolvable lattices.
\end{proposition}


  


\begin{figure}
\begin{center}
\includegraphics[width=8cm]{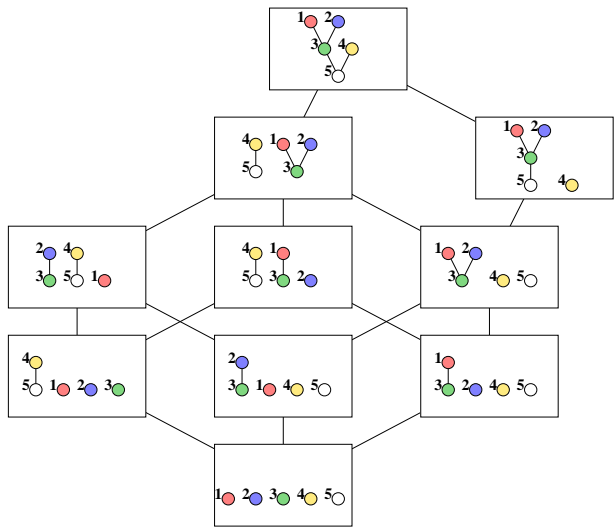} 
\label{interval}
\caption{An interval in the poset $\Pi_{\NAP}(\{1,2,3,4,5\})$.}
\end{center}
\end{figure}

A poset $P$ is called \textbf{totally semi-modular} if for all $x,y \in P$, if
there is $z$ that is covered by both $x$ and $y$, then there is $w$
which covers both $x$ and $y$.

\begin{proposition}
  \label{totalSM}
  The intervals in $\Pi_{\NAP}(n)$ are totally semi-modular lattices.
\end{proposition}

\begin{proof}
  By relation (\ref{intervalles}) and Prop. \ref{subtree}, it is
  enough to prove the proposition for an interval of the form
  $[t,\one]_s$. Let $x,y$ be two sub-rooted trees of $t$. Assume $x$
  and $y$ cover a sub-rooted tree $z$. Hence $x$ is obtained from $z$
  by removing a leaf $l_x$ and $y$ is obtained from $z$ by removing a
  leaf $l_y\not=l_x$. The sub-rooted tree $w$ obtained from $z$ by
  removing the leaves $l_x$ and $l_y$ covers both $x$ and $y$.
\end{proof}


\begin{corollary}
  The posets $\Pi_{\NAP}(n)$ are Cohen-Macaulay.
\end{corollary}
\begin{proof}
  This follows from shellability, hence either from Prop.
  \ref{totalSM}, as total semi-modularity implies CL-shellability
  \cite{BjornerWachs}, or from Prop. \ref{supersolvable}.
\end{proof}

\begin{corollary}
  The operad $\NAP$ is Koszul.
\end{corollary}
\begin{proof}
  This follows from Vallette's criterion Prop. \ref{val_critere} and
  the previous corollary.
\end{proof}

\begin{proposition}
  In the $\NAP$ case, coinvariants are the same as isomorphism classes of
  posets between rooted trees.
\end{proposition}

\begin{proof}
  Coinvariants are given by unlabeled rooted trees. It is clear that if two rooted trees
  have the same underlying unlabeled rooted trees then their
  associated posets are isomorphic. Conversely, let $[\zero,t]$ and
  $[\zero,t']$ be two isomorphic posets. Then $[t,\one]_s$ and
  $[t',\one]_s$ are isomorphic. This isomorphism induces a
  bijection between the labelling of the vertices of the two rooted trees and
  proves that $t$ and $t'$ have the same underlying unlabeled rooted tree.
\end{proof}

This does not work for forests. The Hopf algebra $H_{\NAP}$ is not
free on the coinvariants, as there are relations, given by the
following proposition.

\begin{proposition}
  \label{further_deco}
  Let $t=B(r,t_1,\ldots,t_k)$ be a rooted tree. The poset $[\zero,t]$
  is isomorphic to the product over $j\in\{1,\ldots,k\}$ of the posets
  $[\zero,B(r,t_j)]$.
\end{proposition}
\begin{proof}
By Prop. \ref{subtree} we prove the equivalent result for the
interval  $[t,\one]_s$. Any sub-rooted tree $u$ of $t$ writes
$u=B(r,u_1,\ldots,u_k)$, where $u_i$ is a sub-rooted tree of $t_i$ or
may be the empty tree. The isomorphism sends $u$ to the product of the
rooted trees $B(r,u_j)$.
\end{proof}


\subsection{Isomorphism between  $H_{\NAP}$ and the Hopf algebra of
  Connes and Kreimer}

In \cite{ConKre98}, Connes and Kreimer build a commutative Hopf
algebra $\mathcal H_R$ , polynomial on unlabeled rooted trees. We prove in this section that $H_{\NAP}$ is
isomorphic to this Hopf algebra.

\begin{proposition}\label{libre}
  The Hopf algebra $H_{\NAP}$ is a free commutative algebra on the
  unlabeled rooted trees of root-valence 1.
\end{proposition}

\begin{proof}
  According to \cite[Theorem 6.4]{schmitt}, the Hopf algebra $H_\PP$
  is a free commutative algebra on its set of indecomposable elements.
  But each poset $[\zero,t]$ with $t$ of root-valence $1$ cannot be
  written as a product (because there is only one element covered by
  $\one$), hence is indecomposable. Conversely, any other interval
  decomposes as a product of such intervals by Prop.
  \ref{further_deco}.
\end{proof}

\begin{lemma}
  \label{basisH_NAP}
  The elements $\FF_{[t]}$, where $t$ runs over the set of rooted trees, form a
  basis of the vector space $H_{\NAP}$.
\end{lemma}
\begin{proof}
  As $H_{\NAP}$ is a free algebra on the elements $\FF_{[x]}$ where
  $x$ is an unlabeled rooted tree of root-valence one by
  Prop. \ref{libre}, a vector space basis is
  given by products $\FF_{[t_1]} \dots \FF_{[t_k]}$, with  $t_i=B(r,t'_i)$. By Prop.
  \ref{further_deco}, there exists a unique rooted tree
  $t=B(r,t'_1,\dots,t'_k)$ such that $\prod_i[\zero, t_i]\simeq
  [\zero,t]$. This gives a bijection between forests of unlabeled
  rooted trees of root-valence one and unlabeled rooted trees. Therefore the
  elements $\FF_{[t]}$ where $[t]$ are unlabeled rooted trees form a basis.
\end{proof}

The Hopf algebra of Connes et Kreimer $\mathcal H_R$ is the free
commutative algebra on unlabeled rooted trees with the following coproduct
\begin{equation*}
  \Delta(t)=1\otimes t +t\otimes 1+\sum_c P^c(t)\otimes R^c(t),
\end{equation*}
where $c$ stands for all the admissible cuts, $P^c(t)$ is a forest
and $R^c(t)$ is a rooted tree defined from an admissible cut $c$ (see
\cite{ConKre98}). The coproduct has an alternative definition given by
induction
\begin{equation*}
  \Delta(B^+(t_1,\ldots,t_k))=B^+(t_1,\ldots,t_k)\otimes 1+(\id\otimes
B^+)(\Delta(t_1\ldots t_k)),
\end{equation*}
where $B^+(t_1,\ldots,t_k)=B(r,t_1,\ldots,t_k)$.
This means that the linear map $B^+:\mathcal H_R\rightarrow \mathcal
H_R$ is a $1$-cocycle in the complex computing the Hochschild cohomology
of the coalgebra $\mathcal H_R$. Indeed, Connes and Kreimer prove that
$\mathcal H_R$ is a solution to a universal problem in Hochschild
cohomology.

\begin{theorem}[\cite{ConKre98}]
  The pair $(\mathcal H_R,B^+)$ is universal among commutative Hopf
  algebras $(\mathcal H,L)$ satisfying
  \begin{equation}\label{cocycle}
    \Delta(L(x))=L(x)\otimes 1+(\id\otimes
    L)(\Delta(x)), \forall x\in H.
  \end{equation}
  More precisely, given such a Hopf algebra there exists a unique
  morphism of Hopf algebras $\phi: \mathcal H_R\rightarrow \mathcal H$
  such that $L\circ\phi=\phi\circ B^+$.
\end{theorem}

As a consequence of the universal property, the Hopf algebra of Connes and
Kreimer is unique up to isomorphism. We use this criteria to prove
that $H_{\NAP}$ is isomorphic to $\mathcal H_R$.

\begin{theorem}
  \label{iso_HNAP_HR}
  The Hopf algebra $H_{\NAP}$ is isomorphic to the Hopf algebra
  $\mathcal H_R$ of Connes and Kreimer. The unique isomorphism
  compatible with the universal property sends
  $\FF_{[B(r,t_1,\ldots,t_k)]}$ to the forest $t_1\ldots t_k$.
\end{theorem}

\begin{proof}
  Let us define a $1$-cocycle $L_{\NAP}$ on $H_{\NAP}$. 

  By lemma \ref{basisH_NAP}, it is enough to define $L_{\NAP}$ as
  \begin{equation*}
    L_{\NAP}(\FF_{[t]})=\FF_{[B(r,t)]},
  \end{equation*}
  for each rooted tree $t$.

  Let us prove that $L_{\NAP}$ satisfies equation (\ref{cocycle}).  For any
  rooted tree $u$, let $\psi_u$ be the isomorphism from $[\zero,u]$ to
  $[u,\one]_s$. The coproduct is then given by
  \begin{equation*}
    \Delta(\FF_{[u]})=\sum_{x\leq u} \FF_{[x]}\otimes \FF_{[\psi_u(x)]},
  \end{equation*}
  where $x$ is a forest of rooted trees and $\psi_u(x)$ is a sub-rooted tree of $u$.
  Hence
  \begin{multline*}
    \Delta(\FF_{[B(r,t)]})=\sum_{x\leq B(r,t)} \FF_{[x]}\otimes \FF_{[\psi_{B(r,t)}(x)]} \\
    =\FF_{[B(r,t)]}\otimes 1+\sum_{x\leq t} \FF_{[x]}\otimes \FF_{[\psi_{B(r,t)}(\tilde x)]}.
  \end{multline*}
  Indeed, since $t$ is the unique rooted tree covered by $B(r,t)$, any
  $\tilde x<B(r,t)$ is the forest obtained from a forest of rooted trees $x\leq t$ by
  adding the rooted tree with the single vertex $r$. As a consequence
  $\psi_{B(r,t)}(\tilde x)=B(r,\psi_t(x))$. Hence $L$ satisfies
  equation (\ref{cocycle}).

  \medskip

  Let $(\mathcal H,L)$ be a commutative Hopf algebra satisfying
  relation (\ref{cocycle}). In order to build a morphism of Hopf
  algebras $\rho: H_{\NAP}\rightarrow \mathcal H$, it is
  enough to give its values on the rooted trees of root-valence $1$. But such
  a generator can be written $\FF_{[B(r,t)]}=L_{\NAP}(\FF_{[t]})$. Hence
  we define $\rho(\FF_{[B(r,\emptyset)]})=1$ since the rooted tree with single
  vertex $r$ is the unit and by induction
  $\rho(\FF_{[B(r,t)]})=L(\rho(\FF_{[t]}))$ where $\FF_{[t]}$ is a product
  of generators of degree less than $B(r,t)$.  It is straightforward
  to check that $\rho$ is a morphism of Hopf algebras such that
  $L\circ\rho=\rho\circ L_{\NAP}$.

  As a consequence, $H_{\NAP}$ is isomorphic to $\mathcal H_R$, and
  the isomorphism goes as follows: $\FF_{[B(r,t_1,\ldots,t_k)]}$ is sent
  to the forest $t_1\ldots t_k$.
\end{proof}

Let us give examples for the coproduct $\Delta$ in the incidence
Hopf algebra $H_{\NAP}$:
\begin{equation*}
  \Delta \Farb{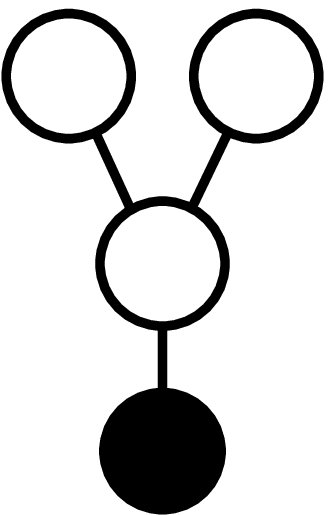}= 1 \otimes \Farb{1200}+2\,\Farb{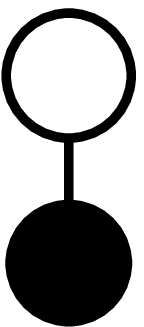}\otimes \Farb{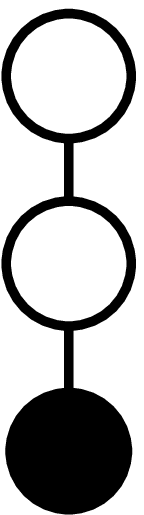}+\Farb{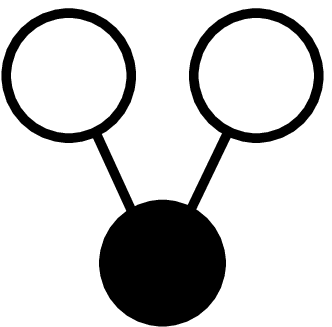}\otimes
  \Farb{10}+\Farb{1200}\otimes 1,
\end{equation*}
\begin{equation*}
  \Delta \Farb{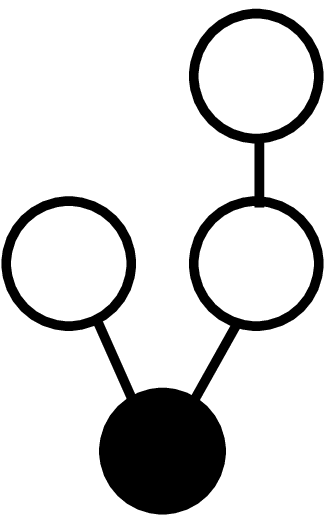}= 1 \otimes \Farb{2100}+\Farb{10}\otimes (\Farb{200}+ \Farb{110})+(\Farb{10}^2+\Farb{110})\otimes \Farb{10}+\Farb{2100}\otimes 1,
\end{equation*}
and
\begin{equation*}
  \Delta \Farb{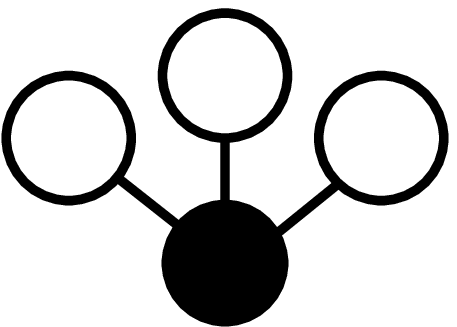}= 1 \otimes \Farb{3000}+3\,\Farb{10}\otimes
  \Farb{200}+3\,\Farb{200}\otimes \Farb{10}+\Farb{3000}\otimes 1,
\end{equation*}
where we have used that
${\FF_{\left[\includegraphics[height=2mm]{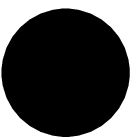}\right]}}$ is
the unit $1$. The last example also follows from the equality (see
Prop. \ref{further_deco})
\begin{equation*}
  \Farb{3000}=\Farb{10}^3.
\end{equation*}

The similar coproducts in the Hopf algebra $\QQ[G_{\NAP}]$ are
\begin{equation*}
  \Delta \Garb{1200}= 1 \otimes \Garb{1200}+\Garb{10}\otimes \Garb{110}+\Garb{200}\otimes
  \Garb{10}+\Garb{1200}\otimes 1,
\end{equation*}
\begin{equation*}
  \Delta \Garb{2100}= 1 \otimes \Garb{2100}+\Garb{10}\otimes (2\, \Garb{200}+ \Garb{110})+(\Garb{10}^2+\Garb{110})\otimes \Garb{10}+\Garb{2100}\otimes 1,
\end{equation*}
and
\begin{equation*}
  \Delta \Garb{3000}= 1 \otimes \Garb{3000}+\Garb{10}\otimes
  \Garb{200}+\Garb{200}\otimes \Garb{10}+\Garb{3000}\otimes 1,
\end{equation*}
where we have also used that ${\GG_{\includegraphics[height=2mm]{a0small.eps}}}$ is the unit $1$.

\subsection{Examples of elements of the  group $G_{\NAP}$}

The group $G_{\NAP}$ can be considered as a group of formal power
series indexed by the set of unlabeled rooted trees. In this section,
we give a criterion for an element of $G_{\NAP}$ to be in $\spec
H_{\NAP}$, then describe several examples of elements of $G_{\NAP}$
and compute their inverses.

Let us first describe explicitly the image of $\spec H_{\NAP}$ in
$G_{\NAP}$.

\begin{lemma}
  \label{decrire_spec}
  A series $f=\sum_{t} \GG_t(f) t$ in $G_{\NAP}$ is in the subgroup $\spec
  H_{\NAP}$ if and only if for each tree $t=B(r,t_1,\dots,t_k)$, one
  has
  \begin{equation*}
    \aut(B(r,t_1,\dots,t_k)) \GG_{B(r,t_1,\dots,t_k)}(f)
    =\prod_{i=1}^k \aut(B(r,t_i))\GG_{B(r,t_i)}(f).
  \end{equation*}
\end{lemma}

\begin{proof}
  Indeed, if the series $f$ is the image of an element $f'$ in $\spec
  H_{\NAP}$, then one has $\GG_t(f)=\FF_{[t]}(f')/\aut(t)$, and the
  multiplicative behavior follows from Prop. \ref{further_deco}.
  Conversely, if the multiplicativity property holds, one can build a
  unique element $f'$ in $\spec H_{\NAP}$ that maps to $f$.
\end{proof}

The first example is in fact an element of the subgroup
$\spec H_{\NAP}$ and we can therefore deduce its inverse by first computing
the Möbius numbers of the maximal intervals in $\Pi_{\NAP}(n)$.

\begin{proposition}
  \label{calcul_mobius}
  Let $t$ be a rooted tree. If $t$ is a corolla with $n+1$ vertices,
  then $\mu(\zero,t)=(-1)^{n}$. If not, then $\mu(\zero,t)=0$.
\end{proposition}

\begin{proof}
  We compute the Möbius number of the poset $[t,\one]_s$.  If $t$
  is the rooted tree $t_2$ with only two vertices then the Möbius number
  of the interval $[t,\one]_s$ is clearly $-1$. Hence the Prop.
  \ref{further_deco} yields the result for the corollas. If the
  valence of the root of $t$ is one then $\one$ covers a unique rooted tree
  which is $t_2$. If $t$ has at least $3$ vertices, then $t_2$ is
  different from $t$ and the Möbius number of the interval
  $[t,1]_s$ is $0$. If the valence of the root of $t$ is greater than
  2 and $t$ is not a corolla, then in the decomposition of
  $t=B(r,t_1,\ldots,t_k)$, there exists $t_i$ having at least two
  vertices. The rooted tree $B(r,t_i)$ has root-valence 1 and has at least $3$
  vertices so its Möbius number is 0. We conclude with Prop. \ref{further_deco}.
\end{proof}

We now deduce from this computation an identity in the group
$G_{\NAP}$.

Consider the series where each rooted tree has weight the inverse of
the order of its automorphism group:
\begin{equation*}
\tt{Z}=\arb{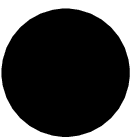}+ \arb{10}+\frac{1}{2}\,\arb{200}+ \arb{110}+\frac{1}{6}\,\arb{3000}+\frac{1}{2}\,\arb{1200}+ \arb{2100}+ \arb{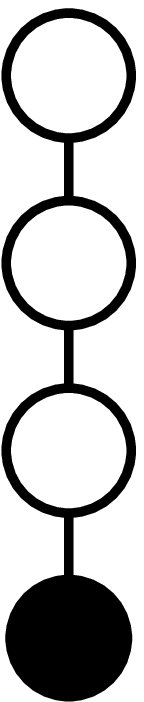}+\dots
\end{equation*}

By Lemma \ref{decrire_spec}, this belongs to the image of $\spec
H_{\NAP}$ and should be called the Zeta function, following the
standard notation in algebraic combinatorics of posets
\cite{rota,stanley1}.

By the general result Prop. \ref{zeta_inverse}, its inverse in the group
$\spec H_{\NAP}$ is known to be the generating series for Möbius
numbers.

Hence by the computation of Möbius numbers done in Prop.
\ref{calcul_mobius} and the inclusion of $\spec H_{\NAP}$ in
$G_{\NAP}$ obtained in Th. \ref{main_thm}, the inverse of $\tt{Z}$
in the group $G_{\NAP}$ is the similar sum $\tt{M}$ restricted on
corollas and with additional signs:
\begin{equation*}
\tt{M}=\arb{0}- \arb{10}+\frac{1}{2}\arb{200}-\frac{1}{6} \arb{3000}+\frac{1}{24}\arb{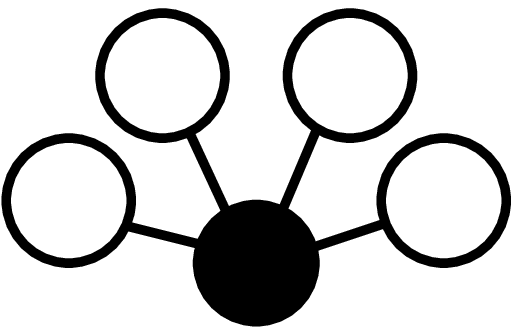}+\dots
\end{equation*}

\medskip

We now give some other examples of elements of $G_{\NAP}$.

Let us introduce the sum of all corollas in $G_{\NAP}$:
\begin{equation*}
\mathtt{C}=\arb{0}+ \arb{10}+\arb{200}+ \arb{3000}+\arb{40000}+\dots
\end{equation*}
and the alternating sum of linear trees:
\begin{equation*}
\mathtt{L}=\arb{0}- \arb{10}+\arb{110}- \arb{1110}+\arb{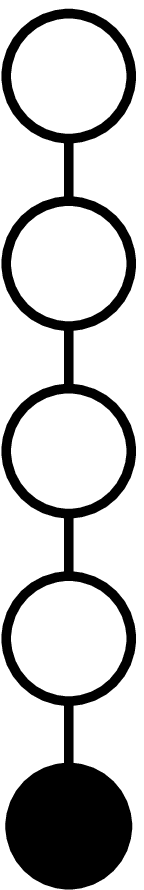}-\dots
\end{equation*}

The series $\mathtt{C}$ satisfies the simple functional equation
\begin{equation}
  \label{eq_fonct_C}
  \mathtt{C}=\arb{0}+\mathtt{C} \tri \arb{0},
\end{equation}
where $\tri$ is the $\NAP$ product on rooted trees.

\begin{theorem}
In the group $G_{\NAP}$, one has $\mathtt{C}=\mathtt{L}^{-1}$.
\end{theorem}
\begin{proof}
  From the functional equation (\ref{eq_fonct_C}) for $\mathtt{C}$,
  one gets
  \begin{equation*}
    \arb{0}=\mathtt{C}^{-1}+\arb{0} \tri \mathtt{C}^{-1},
  \end{equation*} 
  by product by $\mathtt{C}^{-1}$ on the right, since one has in
  $G_{\NAP}$ the relation $(\mathtt{C}\tri \mathtt{D})\times
  \mathtt{E}=(\mathtt{C}\times \mathtt{E})\tri (\mathtt{D}\times \mathtt{E})$. 
But the unique
  solution to this equation is easily seen to be $\mathtt{L}$.
\end{proof}

One can see from Lemma \ref{decrire_spec} that the series $\mathtt{C}$
and $\mathtt{L}$ do not belong to the subgroup $\spec H_{\NAP}$, as
the coefficients $\aut(t) \GG_t$ do not have the necessary
multiplicativity property. For instance, the coefficients of corollas
vanish in $\mathtt{L}$, but the coefficient of the tree $\arb{10}$
does not.

\subsection{Morphisms from $G_{\NAP}$ to usual power series}

There are two morphisms from the group $G_{\NAP}$ to the
\textbf{multiplicative} group of formal power series in one variable
$x$. Either one can project on corollas:
\begin{equation*}
  \sum_{t} \GG_{t}\, t \mapsto \sum_{n \geq 0} \GG_{c_n} x^n,
\end{equation*}
where $c_n$ is the corolla with $n$ leaves, or project on linear trees:
\begin{equation*}
  \sum_{t} \GG_{t}\, t \mapsto \sum_{n \geq 0} \GG_{\ell_n} x^n,
\end{equation*}
where $\ell_{n}$ is the linear tree with $n+1$ vertices. 

Recall that the Hopf algebra of functions on the multiplicative group
of formal power series 
\begin{equation*}
  1+\sum_{n \geq 1} M_n x^n
\end{equation*}
is the free commutative algebra with one generator $M_n$ in each
degree $n \geq 1$ and coproduct
\begin{equation*}
  \Delta M_n =\sum_{i=0}^n M_i \otimes M_{n-i},
\end{equation*}
with the convention that $M_0=1$.

It is indeed easy to check that corollas and linear trees are closed
on the coproduct and that the coproduct is the same as in the
multiplicative group.

In the case of linear trees, the induced morphism from $\spec
H_{\NAP}$ to the multiplicative group of formal power series is again
a projection. This means that it defines a Hopf subalgebra of
$H_{\NAP}$, corresponding to the ladder Hopf subalgebra in $\mathcal
H_R$.

The reader may want to check that the image of the inverse is the
inverse of the image, in the examples $\mathtt{C}$, $\mathtt{L}$,
$\mathtt{Z}$ and $\mathtt{M}$ given above.

\medskip

There is a morphism from the group $G_{\NAP}$ to the group of formal
power series in one variable $x$ for \textbf{composition} given by the
sum of the coefficients of all rooted trees of same degree:
\begin{equation*}
  \sum_{t} \GG_{t}\, t \mapsto \sum_{n \geq 1} \left(\sum_{t \in \NAP(n)_{\sym_n}} \GG_t \right) \, x^{n}.
\end{equation*}
This comes from the morphism of operads from $\NAP$ to the Commutative operad
$\com$ which sends every element of $\NAP(n)$ to the unique
element of $\com(n)$.

One can easily check that the images of the series $\mathtt{C}$ and
$\mathtt{L}$ are inverses for composition. For the images of the
series $\mathtt{Z}$ and $\mathtt{M}$, this is less obvious. This
implies that the image of $\mathtt{Z}$ is
\begin{equation*}
  \sum_{n \geq 1} n^{n-1} \frac{x^n}{n!},
\end{equation*}
which is the functional inverse of $x \exp(-x)$. This is related to
the Lambert W function \cite{lambertW}, the inverse of $x \exp(x)$.

It is easy to see that the induced morphism from $\spec H_{\NAP}$ to
the composition group of formal power series is again a projection,
hence defines a Hopf subalgebra of $H_{\NAP}$. This Hopf algebra is
isomorphic to the Faà di Bruno Hopf algebra as pointed out in \S \ref{S-maintheorem}. 
The generators of the
subalgebra $\QQ[G_{\com}]$ in $\QQ[G_{\NAP}]$ are
\begin{equation*}
  \sum_{t \in \NAP(n)_{\sym_n}} \GG_t,
\end{equation*}
for $n \geq 2$. Hence the generators of the subalgebra in $H_{\NAP}$
are
\begin{equation*}
  \sum_{t \in \NAP(n)_{\sym_n}} \frac{1}{\aut(t)}\FF_{[t]}, 
\end{equation*}
for $n \geq 2$.

Let us give explicitly the first generators of this Hopf subalgebra
of $H_{\NAP}$:
\begin{align*}
 \Gamma_1= &\Farb{10},\\
 \Gamma_2= &\Farb{110}+\frac{1}{2} \Farb{200}=\Farb{110}+\frac{1}{2} \Farb{10}^2,\\
 \Gamma_3= &\Farb{1110}+\frac{1}{2}\Farb{1200}+\Farb{2100}+\frac{1}{6} \Farb{3000}=\Farb{1110}+\frac{1}{2}\Farb{1200}+\Farb{110}\Farb{10}+\frac{1}{6} \Farb{10}^3,
\end{align*}
where we have used the multiplicative property of the $\FF$ functions
to get from sums over all rooted trees to polynomials in rooted trees
of root-valence $1$.

By mapping these elements to the Connes-Kreimer algebra $\mathcal H_R$
by the isomorphism of Theorem \ref{iso_HNAP_HR}, one can see that this
Hopf subalgebra is different from the Connes-Moscovici subalgebra and
also from the other Hopf subalgebras of $\mathcal H_R$ introduced
recently by Foissy \cite{foissy}.







\bibliographystyle{plain}
\bibliography{naposet}

\parbox{.5\linewidth}{
  Muriel Livernet\\
  Laboratoire Analyse, \\Géométrie et Applications,\\
  UMR 7539 du CNRS,\\
  Institut Galilée,\\
  Université Paris Nord\\
  Avenue Jean Baptiste Clément,\\
  93430 VILLETANEUSE.\\
  livernet@math.univ-paris13.fr
}
\parbox{.5\linewidth}{
  Fr{\'e}d{\'e}ric Chapoton\\
  Université de Lyon ;\\
  Université Lyon 1 ;\\
  Institut Camille Jordan CNRS UMR 5208 ;\\
  43, boulevard du 11 novembre 1918,\\
  F-69622 Villeurbanne Cedex.\\
  fchapoton@math.univ-lyon1.fr
}

\end{document}